\documentclass{amsart}
\usepackage{pstricks}

\usepackage{mathtools}

\newtheorem{theorem}{Theorem}

\newcommand{\C}{\mathbb{C}}
\newcommand{\R}{\mathbb{R}}
\renewcommand{\Re}{\operatorname{Re}}
\renewcommand{\Im}{\operatorname{Im}}

\newcommand{\wrt}{\,}
\newcommand{\fn}{\,}

\newcommand{\labs}{\left|}
\newcommand{\rabs}{\right|}
\newcommand{\biglabs}{\bigl|}
\newcommand{\bigrabs}{\bigr|}

\newcommand{\lset}{\left\{}
\newcommand{\rset}{\right\}}
\newcommand{\lpar}{\left(}
\newcommand{\rpar}{\right)}

\newcommand{\biglpar}{\bigl(}
\newcommand{\bigrpar}{\bigr)}

\newcommand{\sgn}{\operatorname{sgn}}
\newcommand{\fudge}{{\frac{1}{\sqrt{2\pi}}}}

\title{On Beurling's Uncertainty Principle}
\author{Xin GAO}
\thanks{This article is part of my Ph.\ D.\ at the University of New South Wales, supervised by Professor Michael Cowling and supported by an Australian Postgraduate Award.}
\address{School of Mathematics and Statistics\\ University of New South Wales\\ UNSW Sydney 2052\\ Australia}

\begin{document}

\begin{abstract}
We generalise a result of Hedenmalm to show that if a function $f$ on $\R$ is such that $\iint_{\R^2} \biglabs f(x) \, \hat f(y) \bigrabs \,e^{\lambda \labs xy \rabs} \,dx\,dy = O( (1-\lambda)^{-N} )$ as $\lambda \to 1-$, then $f$ is the product of a polynomial and a gaussian.
\end{abstract} 
\maketitle
\section{Introduction}

For $a \in \R^+$, define the gaussian $\gamma_a$  by 
\[
\gamma_a(x) = e^{-a x^2/2} 
\qquad \forall x \in \R.
\]
For $f \in L^1(\R)$, define the Fourier transform of $f$ by
\[
\hat f(y) = \fudge \int_\R f(x) \fn e^{- i  x y} \wrt dx
\qquad \forall y \in \R,
\]
and the Mellin transform of $f$ by
\[
M^k_f(z) = \fudge \int_\R f(x) \sgn^k(x) \labs x \rabs^{z - 1/2} \wrt dx 
\]
for $z \in \C$ where this is defined.
Here $k$ takes the values $0$ and $1$.

In 1933, Hardy  \cite{Hardy} proved the following uncertainty principle: if
\[
\labs f \rabs \leq \gamma_a 
\qquad\text{and}\qquad
\biglabs \hat f \bigrabs \leq \gamma_b ,
\]
where $ab = 1$, then $f = c \, \gamma_a$ for some constant $c$.
As a consequence, we obtain the weak uncertainty principle: if $\labs f \rabs \leq \gamma_a$ and $| \hat f | \leq \gamma_b$, where $ab>1$, then $f(x) = 0$ (almost everywhere).

At an unknown time (see Hormander \cite{Hormander} for more information), Beurling proved another uncertainty principle: if $f$ is integrable and
\[
\iint_{\R^2} \biglabs f(x) \, \hat f(y) \bigrabs \,e^{\labs xy \rabs} \,dx\,dy < \infty ,
\]
then $f = 0$.
This result implies the weak form of Hardy's uncertainty principle immediately, and the strong form with a bit more work.
See for instance the last part of \cite{CEKPV} for more on this.

In 2003, Bonami, Demange and Jaming \cite{BDJ} (see also Bonami and Demange \cite{BD}) strengthened both principles, and showed that if
\begin{equation}\label{eq:BDJ-hypo}
\iint_{\R^2} \frac{ \biglabs f(x) \, \hat f(y) \bigrabs}{(1 + \labs x \rabs + \labs y \rabs)^N} \,e^{\labs xy \rabs} \,dx\,dy < \infty ,
\end{equation}
where $N > 0$, then there exist $a\in \R^+$ and a polynomial $p$ such that 
\[
f = p \, \gamma_a.
\]
This implies both Hardy's and Beurling's results.

In 2011, Hedenmalm \cite{Hedenmalm}  came up with a new version of the uncertainty principle, as follows: if 
\[
\iint_{\R^2} \biglabs f(x) \, \hat f(y) \bigrabs \,e^{\lambda \labs xy \rabs} \,dx\,dy = O( (1-\lambda)^{-1} ) 
\qquad\text{as $\lambda \to 1-$},
\]
then there is a constant $c$ such that $| M^k_f(z) | = c | M^k_{\gamma_1}(z)| $ when $z \in i \R$.
Note that $| M^k_{\gamma_a}(z) | = c | M^k_{\gamma_1}(z)|$ when  $z \in i \R$.
This theorem implies Beurling's result, but not that of Bonami \textit{et al.} nor that of Hardy, or at least not immediately.

In this paper, we further extend these uncertainty principles, as follows.
\begin{theorem}
If 
\begin{equation}\label{eq:hypo}
\iint_{\R^2} \biglabs f(x) \, \hat f(y) \bigrabs \,e^{\lambda \labs xy \rabs} \,dx\,dy = O( (1-\lambda)^{-N} ) 
\qquad\text{as $\lambda \to 1-$},
\end{equation}
where $N > 0$, then there exist $a \in \R^+$ and a polynomial $p$ such that 
\[
f = p \, \gamma_a.
\]
\end{theorem}
This result implies Beurling's, Hedenmalm's and Hardy's results, and, the proof also yields that of Bonami \textit{et al}.
We plan to discuss the extension to $\R^d$ in a future article.

\section{The proof}
There are a number of steps in the proof. 
In Steps 1 and 2, we observe that it suffices to suppose that $f$ is odd or even, and real-valued, and show that $f$ vanishes at infinity faster than any exponential.
In Steps 3 and 4, we introduce an auxilliary function $F$, similar to that introduced by Hedenmalm \cite{Hedenmalm}, and then find its Mellin transform.
Step 5 finds the Mellin transform of $f$ and Step 6 finds $f$.

\subsection{Step 1}
Observe that if $f $ may be written as a sum, that is $p \, \gamma_a + q \, \gamma_b$, where $p$ and $q$ are polynomials and $0 < a \leq b$, then $\hat f = \tilde p \, \gamma_{1/a} + \tilde q \, \gamma_{1/b}$, where $\tilde p$ and $\tilde q$ are also polynomials and $0 < 1/b \leq 1/a$.
Therefore there exists a positive constant $R$ such that $\biglabs f(x)\bigrabs \geq c \, \gamma_a$ when $x \geq R$ and $\biglabs \hat f(y)\bigrabs \geq d \, \gamma_{1/b}$ when $y \geq R$.
Further, the first quadrant is the disjoint union of three sets; the first where $0 \leq x < R$ and $y \geq 0$; the second where $x \geq R$ and $0 \leq y < R$, and the third where $x \geq R$ and $y \geq R$.
When $0 \leq \lambda \leq 1$, it is clear that
\[
 \int_0^\infty \int_0^R \gamma_a(x) \, \gamma_{1/b}(y) \, e^{\lambda xy} \,dx \,dy 
 \leq R \int_{0}^\infty \gamma_{1/b}(y) \, e^{Ry}  \,dy 
 < \infty
\]
and
\[
 \int_{0}^{R} \int_R^\infty \gamma_a(x) \, \gamma_{1/b}(y) \, e^{\lambda xy} \,dx \,dy 
\leq \int_{R}^{\infty} \gamma_a(x) \, e^{Rx} \,dx 
< \infty ,
\]
and so
\[
 \int_{R}^{\infty} \int_R^\infty \gamma_a(x) \, \gamma_{1/b}(y) \, e^{\lambda xy} \,dx \,dy 
 < \infty
\]
if and only if 
\[
  \int_{0}^{\infty} \int_0^\infty \gamma_a(x) \, \gamma_{1/b}(y) \, e^{\lambda xy} \,dx \,dy 
 < \infty.
\]
Moreover, 
\[
\begin{aligned}
 \int_{0}^{\infty} \int_0^\infty \gamma_a(x) \, \gamma_{1/b}(y) \, e^{\lambda xy} \,dx \,dy 
& =   \int_{0}^{\infty} \int_0^\infty e^{-a(x- \lambda y/a)^2/2} \, e^{\lambda^2 y^2/2a -y^2/2b} \,dx \,dy \\
& =   \int_{0}^{\infty} \int_0^\infty \gamma_a(x- \lambda y/a) \, \gamma_{c}(y) \,dx \,dy \\
& \geq \int_{0}^{\infty} \int_0^\infty \gamma_a(u) \, \gamma_{c}(y) \,du \,dy ,
 \end{aligned}
\]
where $c = 1/b-\lambda^2/a$; this integral is finite for all $\lambda \in [0,1)$ if and only if $a/b \geq 1$.
So \eqref{eq:hypo} implies that $a = b$.
Similar arguments work if $f$ is a finite sum of terms $p \, \gamma_a$ and also show that $a$ cannot be complex.
A more refined version of this argument allows us to control the degree of the polynomial $p$ by the number $N$ in \eqref{eq:hypo}.

The point is that if $f$ satisifies the hypothesis \eqref{eq:hypo}, then so do its odd and even parts and its real and imaginary parts.
Consequently, there is no loss of generality in supposing that $f$ is real-valued, and either odd or even.

\subsection{Step 2}
If $f \neq 0$ and
\[
\iint_{\R^2} \biglabs f(x) \, \hat f(y) \bigrabs \,e^{\lambda \labs xy \rabs} \,dx\,dy = O( (1-\lambda)^{-N} ) 
\qquad\text{as $\lambda \to 1-$},
\]
then $f$ and $\hat f$ are integrable, whence both are continuous. 
Thus 
\begin{equation}\label{eq:exp-decay}
\int_\R  \labs f(x) \rabs \,e^{\epsilon \labs x \rabs} \,dx < \infty
\qquad\text{and}\qquad
\int_\R \biglabs \hat f(y) \bigrabs \,e^{\epsilon \labs y \rabs} \,dy < \infty 
\end{equation}
for some $\epsilon \in \R^+$.
Hence $\hat f$ and $f$ are real analytic, which implies that neither $f$ nor $\hat f$ is compactly supported.
We may now conclude that \eqref{eq:exp-decay} is true for all $\epsilon \in \R^+$.

\subsection{Step 3}
We are going to use a variant of the auxillary function defined by Hedenmalm; we define $F:\R \to \C$ as follows:
\begin{equation*}
F(\lambda) = \fudge \int_\R f(x) \, f(\lambda x) \, dx
\qquad\forall \lambda \in \R.
\end{equation*}
(Hedenmalm had a complex conjugate on $f(x)$.)
In this step, we prove that
\[
\begin{aligned}
F(\lambda) 
&= \sum_{j=0}^{2N} \frac{t_j}{(\lambda^2+1)^{j+1/2}} + \frac{u_j \lambda}{(\lambda^2+1)^{j+1/2}} .
\end{aligned}
\]
where $t_j$ and $u_j$ are complex numbers, and $u_0 = 0$.

A change of variables shows that
%
%
%
\[
F\lpar \frac{1}{\lambda}\rpar = \labs \lambda\rabs F(\lambda)
\qquad\forall \lambda \in \R \setminus \{0\}.
\] 
We set $G(\lambda) = (\lambda^2 + 1) ^{1/2} \, F(\lambda)$, and then
\[
\begin{aligned}
G\lpar \frac{1}{\lambda}\rpar 
&= G(\lambda).
\end{aligned}
\]

By Fourier inversion
\[
F(\lambda) = \frac{1}{2\pi} \iint_{\R^2} f(x) \, \hat f(y) \, e^{i\lambda xy} \,dx\,dy
\qquad\forall \lambda \in \R,
\]
so the hypothesis \eqref{eq:hypo} implies that $F$ and hence $G$ extend analytically to the strip $\{ \lambda \in \C : \labs \Im \lambda \rabs < 1 \}$.
Therefore the function $\lambda \mapsto G(1/\lambda)$ extends analytically to $\C$ less two closed circles of radius $1/2$ with centres at $\pm i/2$, and coincides with $G$ on the real axis and hence where both are defined.
Thus $G$ extends to an analytic function on $\C \setminus \{\pm i\}$.

We claim that
\[
\labs G(\lambda) \rabs = O( \labs\lambda^2 + 1\rabs^{- 2N + 1/2} )  \qquad\text{as $\lambda \to \pm i$.}
\]
To prove this, we consider $\lambda$ lying inside the unit ball $B(0,1)$; the reflection property gives the result for $\lambda$ lying outside the ball.
By assumption, if $\lambda \in B(0,1)$, then
\[
\begin{aligned}
\labs G(\lambda) \rabs 
&= \labs\lambda^2 + 1\rabs^{1/2} \labs F(\lambda) \rabs \\
&\leq \frac{\labs\lambda^2 + 1\rabs^{1/2} }{2\pi} \labs \iint_{\R^2}  \ f(x) \, \hat f(y) \, e^{i\lambda xy} \,dx\,dy \rabs\\
&\leq \frac{\labs\lambda^2 + 1\rabs^{1/2} }{2\pi} \iint_{\R^2} \biglabs \ f(x) \, \hat f(y) \bigrabs e^{\labs \Im\lambda xy\rabs} \,dx\,dy \\
&=  \labs\lambda^2 + 1\rabs^{1/2} O( (1 - \labs\Im\lambda\rabs)^{- N} )\\
&= O\biglpar \labs\lambda^2 + 1\rabs^{- 2N + 1/2} \bigrpar.
\end{aligned}
\]

From this inequality, $G$ is meromorphic.
We deduce that
\[
\begin{aligned}
G(\lambda) 
&= \sum_{j=0}^{2N} \frac{c_j}{(\lambda -i)^j} + \frac{d_j}{(\lambda +i)^j}  
  = \sum_{j=0}^{2N} \frac{c_j (\lambda +i)^j + d_j(\lambda -i)^j }{(\lambda ^2+1)^j} \\
&= \sum_{j=0}^{2N} \frac{p_j(\lambda ^2) + \lambda  q_j(\lambda ^2)}{(\lambda ^2+1)^j} 
   = \sum_{j=0}^{2N} \frac{r_j(\lambda ^2+1) + \lambda  s_j(\lambda ^2+1)}{(\lambda ^2+1)^j} \\
&= \sum_{j=0}^{2N} \frac{t_j}{(\lambda ^2+1)^j} + \frac{u_j \lambda }{(\lambda ^2+1)^j} ,
\end{aligned}
\]
where $p_j$, $q_j$, $r_j$ and $s_j$ are polynomials and $c_j$, $d_j$, $t_j$ and $u_j$ are complex numbers.
Thus
\[
\begin{aligned}
F(\lambda) 
&= \sum_{j=0}^{2N} \frac{t_j}{(\lambda^2+1)^{j+1/2}} + \frac{u_j \lambda}{(\lambda^2+1)^{j+1/2}} .
\end{aligned}
\]
Note that $u_0 = 0$ since $F$ vanishes at $\infty$.

\subsection{Step 4}
A routine computation of the Mellin transform shows that
\begin{equation}\label{eq:comp-M0F}
\begin{aligned}
M^0_F (z)
&= \sum_{j=0}^{2N} \fudge t_j  \, \frac{ \Gamma(z/2 + 1/4) \, \Gamma (j - z/2 + 1/4) }{ \Gamma(j + 1/2)}   \\
&= p^0(z) \, \Gamma(z/2 + 1/4) \, \Gamma ( - z/2 + 1/4) ,
\end{aligned}
\end{equation}
and
\begin{equation}\label{eq:comp-M1F}
\begin{aligned}
M^1_F (z)
&= \sum_{j=1}^{2N} \fudge u_j  \, \frac{ \Gamma(z/2 + 3/4) \, \Gamma (j - z/2 - 1/4) }{ \Gamma(j + 1/2)}   \\
&= p^1(z) \, \Gamma(z/2 + 3/4) \, \Gamma ( - z/2 + 3/4) .
\end{aligned}
\end{equation}
Here, $p^0$ and $p^1$ are polynomials.

\subsection{Interlude}
The Mellin version of the Fourier transform of a convolution being the product of the Fourier transforms is that
\begin{equation}\label{eq:MF-Mf}
M^k_F (z) =  M^k_f(z) \, M^k_f(-z)
\qquad\forall z \in i \R .
\end{equation}
We now have information about $M^k_F(z)  = M^k_f(z) \, M^k_f(-z)$.
Hedenmalm defined $F$ with a complex conjugate and extracted information about $| M^k_f(z) |^2$ when $z \in i\R$ in much the same way.
The problem is to turn this into information about $M^k_f(z)$, for this will determine $f$ by an inverse Mellin transformation.

Define 
\begin{equation}\label{eq:def-Theta}
\Theta^k_f(z)  =  \frac{M^k_f(z)}{\Gamma( z/2 + k /2 + 1/4)}   
\end{equation}
where $k = 0, 1$, for $z$ for which this makes sense.
We claim that $\Theta^k_f(z)$ extends to an entire function, and there exist constants $c$ and $d$ such that
\[
\biglabs \Theta^k_f(z) \bigrabs \leq c \,e^{d \labs z\rabs \ln(1+\labs z \rabs)} .
\]
We also claim that
\[
\Theta^k_f(z) \, \Theta^k_f(-z) = p^k(z) 
\qquad\forall z \in \C.
\]
It follows that $\Theta^k_f$ is an entire function with at most $N$ zeros, of exponential type $1$, and hence 
\[
\Theta^k_f(z) = r^k(z) \, e^{a^k z} ,
\]
where $r^k$ is a polynomial of degree at most $N$ and $a^k \in \C$.

This enables us to pin down $f$ as a sum of the form 
\[
s \,\gamma_b + t \, \gamma_c, 
\]
where $s$ and $t$ are polynomials, one even and the other odd.
From Step 1, $b = c$.

\subsection{Step 5}
Recall that
\[
\Theta^k_f(z)  =  \frac{M^k_f(z)}{\Gamma( z/2 + k /2 + 1/4)}   
\]
where $k = 0, 1$, for $z$ for which this makes sense.
We will show that $\Theta^k_f(z)$ extends to an entire function, and there exist constants $C$ and $D$ such that
\begin{equation}\label{eq:order-1}
\biglabs \Theta^k_f(z) \bigrabs \leq C e^{D \labs z\rabs \ln(1 + \labs z \rabs)} .
\end{equation}

When $\Re(z) \geq 0$, this follows from  the estimate 
\[
\begin{aligned}
\biglabs M^k_f(z) \bigrabs &\leq \int_\R \labs f(x) \rabs \labs x \rabs ^{\Re z -1/2} \,dx \\
&\leq C \int_{\R^+} e^{ \labs x \rabs} \labs x \rabs ^{\Re z -1/2} \,dx \\
&= C \,  \Gamma(\Re z + 1/2);
\end{aligned}
\]
the absolute convergence of the integral gives the analytic continuation and the estimate \eqref{eq:order-1} follows from the asymptotic formula for the gamma function.

Observe that, when $\labs \Re z\rabs < 1/2$, 
\[
\begin{aligned}
\Theta^k_{\hat f}(z)
&= \int_\R \hat f(y) \frac{\sgn(y)^k \labs y \rabs^{z - 1/2} }{\Gamma( z/2 + k /2 + 1/4)} \,dy  \\
&= \int_\R  f(x) \lpar \frac{\sgn(\cdot)^k \labs \cdot \rabs^{z -1/2} }{\Gamma( z/2 + k /2 + 1/4)}\rpar\!\!\widehat{\phantom{M}} (x) \,dx \\
&= i^{-k} 2^{z} \int_\R  f(x) \, \frac{\sgn(\cdot)^k \labs x \rabs^{-z -1/2} }{\Gamma( -z/2 + k /2 + 1/4)} \,dx \\
&= i^{-k} 2^{z} \Theta^k_{f}(-z);
\end{aligned}
\]
the last term continues analytically into the left half-plane, and we can estimate $\Theta^k_f(z)$ when $\Re z \leq 0$ by estimating $\Theta^k_{\hat f}(z)$ when $\Re z \geq 0$.

\subsection{Step 6}
We have shown that $\Theta^k_f(z)$ extends to an entire function of order less than 2(see \eqref{eq:order-1}).
Moreover, by \eqref{eq:def-Theta}, \eqref{eq:MF-Mf}, \eqref{eq:comp-M0F}, and \eqref{eq:comp-M1F},
\begin{equation}\label{eq:Theta-Theta}
\begin{aligned}
\Theta^k_f(z) \, \Theta^k_f(-z) 
&=  \frac{M^k_f(z)}{\Gamma( z/2 + k /2 + 1/4)} \,\frac{M^k_f(-z)}{\Gamma( -z/2 + k /2 + 1/4)} \\
&= \frac{M^k_F(z)}{\Gamma( z/2 + k /2 + 1/4) \, \Gamma( -z/2 + k /2 + 1/4)} \\
&= p^k(z) 
\qquad\forall z \in \C,
\end{aligned}
\end{equation}
where $p^k$ is a polynomial thus has finite zeros.
As we remarked earlier, it follows that $\Theta^k_f$ is an entire function of order $1$ with at finitely many zeros, and hence 
\[
\Theta^k_f(z) = r^k(z) \, e^{a^k z} ,
\]
where $r^k$ is a polynomial and $a^k \in \C$, by the Hadamard Factorisation Theorem (see, e.g., \cite[\S8.24]{Titchmarsh}).
Moreover, by Step 1, we may and shall suppose that $f$ is real-valued, and so when $z$ is purely imaginary,
\[
\labs \Theta^k_f(z) \rabs = \labs \frac{ M^k_f(z)}{\Gamma(\frac{1}{4} + \frac{z}{2})}\rabs = \labs \frac{ M^k_f(-z)}{\Gamma(\frac{1}{4} - \frac{z}{2})}\rabs = \labs \Theta^k_f(-z) \rabs ,
\]
whence $| \Theta^k_f(z) |$ grows polynomially in $|z|$ when $\Re(z) =0$, from \eqref{eq:Theta-Theta}.
This implies that $a^k$ is real.

Now we invert the Mellin transformation.
In light of Step 1, it suffices to treat only $f_o$, the odd part of $f$; the argument for the even part is similar.
It is easy to check that if $\delta_a(x) = x \,\gamma_a(x)$, then $M_{\delta_a}^0 = 0$, while
\[
\begin{aligned}
M_{\delta_a}^1(z) 
&= \frac{1 }{\sqrt{2\pi}} \, \lpar \frac{a}{2} \rpar^{-z/2 - 3/4} \, \Gamma(z/2 + 3/4) .
\end{aligned}
\] 
By taking the partial derivatives of both sides with respect to $a$, we see that the Mellin transform of the function $\delta^k_a: x \mapsto x^{k} \, e^{-ax^2/2}$ is equal to $0$ if $k$ is even and $q^k(z) \lpar {a}/{2} \rpar^{-z/2 - 3/4} \, \Gamma(z/2 + 3/4)$, where $q$ is a polynomial of degree exactly $(k-1)/2$ if $k$ is odd.
By taking suitable linear combinations of these formulae, we may show that there is a polynomial $p$ such that $f = p \, \gamma_b$, where $b = \ln(a/2)$.

\section{Remarks}
If we make the assumption of Bonami, Demange and Jaming (see \eqref{eq:BDJ-hypo}), we can still go through the arguments of Step 2 to show that $f$ and $\hat f$ decay faster than any exponential at infinity.
This implies that
\[
\iint_{  \min\lset|x|, |y| \rset \leq 1 }  \biglabs f(x) \,\hat f(y) \bigrabs \,e^{\lambda \labs xy \rabs} \,dx\,dy < \infty
\]
for all $\lambda \in \R^+$.
Moreover, 
\[
\begin{aligned}{}
&\iint_{ \min\lset|x|, |y| \rset \geq 1 }  \biglabs f(x) \,\hat f(y) \bigrabs \,e^{\lambda \labs xy \rabs} \,dx\,dy  \\
&\quad= O((1-\lambda)^N ) \iint_{ \min\lset|x|, |y| \rset \geq 1}  \frac{ \biglabs f(x) \,\hat f(y) \bigrabs } {(1 + \labs x \rabs + \labs y \rabs)^N } \,e^{ \labs xy \rabs} \,dx\,dy  ,
\end{aligned}
\]
since
\[
\max_{ \min\lset|x|, |y| \rset \geq 1} e^{ (\lambda - 1) \labs xy \rabs } (1 + |x| + |y|)^N  = O((1-\lambda)^N ) .
\]
Indeed, on a line segment where $|x| + |y|$ is constant, $e^{ (\lambda - 1) \labs xy \rabs } $ is maximum when $|x|= 1$ or $|y| =1$.

Hence our arguments also imply the result of Bonami, Demange, and Jaming.

\end{document}